\documentclass[leqno,12pt]{article}%
\usepackage{amsmath,amssymb}
\usepackage{dsfont}
\usepackage{amsfonts}%
\usepackage{color}
\usepackage{epsfig}
\usepackage{pgf,pgfarrows,pgfnodes,pgfautomata,pgfheaps}

\newtheorem{Theorem}{\hspace{\parindent}\bf Theorem}[section]

\newtheorem{Proposition}{\hspace{\parindent}\bf Proposition}[section]
\newtheorem{Corollary}{\hspace{\parindent}\bf Corollary}[section]

\newcommand{\qed}{\hfill$\square$\vspace{0.3cm}}

\pgfdeclareimage[height=5.cm]{fig1}{Regiones}

\begin{document}

\title{\textbf{Anisotropic nonlocal diffusion equations with singular forcing}}
\author{by\\Arturo de Pablo, Fernando Quir\'{o}s, and Ana Rodr\'{\i}guez}

\maketitle

\

\begin{abstract}
We prove existence, uniqueness and  regularity of solutions of nonlocal heat equations associated to anisotropic stable diffusion operators. The main features are that the right-hand side has  very few regularity and that the spectral measure can be singular in some directions. The proofs require having good enough estimates for the corresponding heat kernels and their derivatives.
\end{abstract}


\vskip 8cm

\noindent{\makebox[1in]\hrulefill}\newline2010 \textit{Mathematics Subject
Classification.} 35R11, 
35B65, 
35A05. 
\newline\textit{Keywords and phrases.} Non-local diffusion, anisotropic stable operators, well-posedness, regularity, singular forcing, heat kernel estimates for L\'evy processes.

\newpage

\section{Introduction and main results}

\label{sect-introduction} \setcounter{equation}{0}

The aim of this paper is to study existence, uniqueness and regularity of solutions $u$ to a nonlocal parabolic problem with a nonstandard forcing term,
\begin{equation}  \label{eq:main}
\left\{
\begin{array}{ll}
\partial_t u+\mathcal{L}u=\mathcal{L}f\quad & \text{in }  Q:=\mathbb{R}^N\times\mathbb{R}_+,
\\ [4mm]
u(\cdot,0) = u_0\quad &\text{in }\mathbb{R}^N,
\end{array}
\right.
\end{equation}
where $f=f(x,t)$ is bounded and H\"older continuous,  $u_0$ is in some $L^p$ space, $1\le p\le \infty$, and $\mathcal{L}$ is a pseudo-differential operator corresponding to a symmetric stable process of order $\sigma\in(0,2)$. In general the right-hand side $\mathcal{L}f$ is singular since it is  only known to be a distribution in some negative H\"older space. As a consequence,  $u$ is not an energy solution. Thus, though the equation is linear and invariant under translations and scalings, de Giorgi or Moser-like approaches to regularity cannot be applied.

Parabolic equations with singular forcing terms have attracted a lot of attention in recent years; see for instance~\cite{Hairer-2011, Otto-Sauer-Smith-Weber-2018}.  Problem~\eqref{eq:main} has already been considered in~\cite{pqrv3, VPQR} in the special case in which $\mathcal{L}$ is the fractional Laplacian as a tool to obtain higher regularity for nonlinear problems of the form $\partial_t u+(-\Delta)^{\sigma/2}\phi(u)=0$. Hence, we expect that the results in this paper will allow to obtain higher regularity for nonlinear problems associated to more general nonlocal operators.

The nonlocal operator $\mathcal{L}$ is defined, for regular functions which do not grow too much at infinity, by
\begin{equation}
\label{eq:general.operator}
\mathcal{L}u(x) =
\int_{\mathbb{R}^N}\left(u(x)-\frac{u(x +y) + u(x-y)}2\right)
\,d\nu(y),
\end{equation}
where the nonnegative L\'evy measure $\nu$ has the polar decomposition
\begin{equation}
\label{eq:spectralmeasutre}
d\nu(y)=\frac{d(|y|)}{|y|^{N+\sigma}}d\mu\left(\frac y{|y|}\right), \qquad\text{with }0<\sigma<2.
\end{equation}

General operators of the form \eqref{eq:general.operator}--\eqref{eq:spectralmeasutre} arise as the infinitesimal generators of symmetric stable L\'evy processes $X=\{X_t\}_{t\ge0}$, which satisfy
$$
\lambda X_t= X_{\lambda^\sigma t},\qquad \lambda>0,\; t\ge0.
$$
These processes appear  in Physics, Mathematical Finance and Biology, among other applications, and have been the subject of intensive research in the last years from the point of view both of Probability and Analysis; see for instance the survey \cite{Ros-2016} and the references therein.

 The measure $\mu$ on the sphere $\mathbb{S}^{N-1}$, called the \emph{spectral measure}, is assumed to be  finite, $\mu(\mathbb{S}^{N-1})=\Lambda<\infty$,  and to satisfy
the \lq\lq ellipticity'' (non-degeneracy) condition
\begin{equation}\label{eq.ellipcity}
 \inf_{\zeta\in \mathbb{S}^{N-1}} \int_{\mathbb{S}^{N-1}}|\zeta\cdot\theta|^{\sigma}\,d\mu(\theta)\ge\lambda>0.
\end{equation}
That is, we  require that the spectral measure  is not supported in any proper subspace of $\mathbb{R}^N$.
We remark that we do not impose any symmetry to the measure $\mu$, since symmetry of the operator (and thus of the process) comes directly from the way we write it in~\eqref{eq:general.operator} using second differences. If $\mu$ were symmetric the operator would take the more familiar form
$$
\mathcal{L}u(x) =P.V.\int_{\mathbb{R}^{N}}\left(u( x)-u(x + y)\right)\,d\nu(y).
$$
By~\eqref{eq:spectralmeasutre}, the operator $\mathcal{L}$ can be expressed in the form
\begin{equation}
\label{eq:operador en polares}
\mathcal{L}u(x) =\int_{\mathbb{S}^{N-1}}\int_0^\infty
\left(u( x)-\frac{u(x + r\theta) + u(x-r\theta)}2 \right)
\frac{dr}{r^{1+\sigma}}\,d\mu(\theta).
\end{equation}

The spectral measure is allowed to be \emph{anisotropic}. Hence,  we cannot use radial arguments as the ones employed in~\cite{VPQR} to deal with the isotropic case $d\mu(\theta)=d\theta$, for which  the operator reduces to (a multiple of) the well known fractional Laplacian, $\mathcal{L}=(-\Delta)^{\sigma/2}$.   Note, however, that our anisotropic operator is still homogeneous of order $\sigma$, which will turn out to be an important tool in our proofs. We also allow the L\'evy  measure to be singular in some directions. It may even be concentrated on a set of directions of Lebesgue measure zero, as in example~\eqref{eq:sumalaplas}.

The functions $u$ and $f$ in equation \eqref{eq:main} do not possess in general the required regularity to give   $\mathcal{L}u$ and $\mathcal{L}f$ a pointwise meaning. Moreover, $u$ will not even belong to the energy space associated to the operator. We therefore have to work with solutions in a very weak sense; see formulas~\eqref{veryweak0} and~\eqref{veryweak} below. Our first result shows that the problem is well-posed in this very weak formulation.

\begin{Theorem}
\label{th:main}   If $f\in C^\alpha(Q)\cap L^\infty(Q)$ for some $0<\alpha<1$, and $u_0\in
L^p(\mathbb{R}^N)$ for some $1\le p\le\infty$, then  there exists a unique
very weak solution of problem~\eqref{eq:main}, which is  moreover bounded for positive times.
\end{Theorem}

The solution is given explicitly by Duhamel's type formula
\begin{equation}
\label{eq:Duhamel.intro}
u(x,t)=\int_{\mathbb{R}^N}P(x-\overline x,t)u_0(\overline x)\,d\overline x+\int_0^t\int_{\mathbb{R}^N}\mathcal{L}P(x-\overline x,
t-\overline t) f(\overline x,\overline t)\,d\overline xd\overline t,
\end{equation}
where $P$ is the fundamental solution for the homogeneous problem. Note that in the second term the operator is applied to $P$, in contrast with the case in which the right-hand side  is not singular. Since $P$ is smooth for positive times, the first term, corresponding to the initial datum,  is smooth. On the contrary, giving a meaning to the second integral in~\eqref{eq:Duhamel.intro} in some principal value sense, see~\eqref{eq:duhamel1}, requires some effort, since $\mathcal{L}P$ has a nonintegrable singularity at $(x-\bar x,t-\bar t)=(0,0)$. However, thanks to  the homogeneity of the operator, the kernel and its derivatives are known to have a self-similar structure. In particular,  $P(x,t)=t^{-\frac N\sigma}\Phi(xt^{-\frac1\sigma})$, $\mathcal{L}P(x,t)=t^{-\frac N\sigma-1}\mathcal{L}\Phi(xt^{-\frac1\sigma})$ for some positive profile $\Phi$. The singularity of $\mathcal{L}P$ at the origin can then be controlled in terms of the decay of the profile $\Phi$ at infinity, which is next combined with the H\"older regularity of $f$ to obtain integrability at the origin. The same decay estimates account for integrability at space infinity.

In the isotropic case, $\mathcal{L}=(-\Delta)^{\sigma/2}$, we have the pointwise estimates
\begin{equation}\label{eq:first-decay}
\Phi(r)\sim r^{-N-\sigma},\quad |\mathcal{L}\Phi(r)|\le cr^{-N-\sigma}\quad\text{for }r \text{ large},
\end{equation}
which have been known for a long time~\cite{Blumenthal-Getoor}. However, for anisotropic processes the pointwise decay can be much slower in some directions, as observed in~\cite{Pruitt-Taylor} for the  case
\begin{equation}
  \label{eq:sumalaplas}
  \mathcal{L}=\sum_{j=1}^N (-\partial_{x_jx_j})^{\sigma/2};
\end{equation}
see also Section~\ref{sec-suma}.
Nevertheless,  the decay  for $\Phi$  \emph{on average}
\begin{equation}\label{profile-pruitt}
\int_{\mathbb{S}^{N-1}}\Phi(r\theta)\,d\theta\le c r^{-N-\sigma}
\end{equation}
holds for all stable  processes~\cite[Theorem 2]{Pruitt-Taylor}. On the other hand,  for any $k\in\mathbb{N}$ and $\beta\in(0,\sigma)$ there is a function $\Omega_{k,\beta}\in L^1(\mathbb{S}^{N-1})$  such that
\begin{equation}
\label{eq:fractional.derivatives}
|\mathcal{L}^k\Phi(r\theta)|\le \Omega_{k,\beta}(\theta) r^{-N-\beta};
\end{equation}
see~\cite[Theorem 5.1]{Glowacki-1993}. Hence,
\begin{equation}
  \label{eq:decay.time}
  \int_{\mathbb{S}^{N-1}}|\mathcal{L}\Phi(r\theta)|\,d\theta\le c_\beta r^{-N-\beta}\quad \text{for every } \beta<\sigma.
\end{equation}
 The  estimates on average~\eqref{profile-pruitt} and~\eqref{eq:decay.time} suffice to show the well-posedness of our problem.

\medskip

\noindent\emph{Remark. } When $k=0$ an estimate like~\eqref{eq:fractional.derivatives}  with $\beta=\sigma$ is only true if the spectral measure is absolutely continuous and its density belongs to a certain integrability class~\cite{Glowacki-Hebisch}. However, the threshold decay is reached on average; see~\eqref{profile-pruitt}. We expect to have also such limit decay on average for $k\ge 1$.

\medskip

In order to prove that $u$ is H\"older continuous we cannot use de Giorgi or Moser approaches, as done for instance in~\cite{Felsinger-Kassmann-2013, Kassmann-Schwab-2014}, since  the solution does not lie in general in the energy space. Hence we have chosen a different approach, that requires to estimate further derivatives of $\mathcal{L}P$, a subject that has independent interest. The pointwise  estimate~\eqref{eq:fractional.derivatives} provides the required decay on average for $\mathcal{L}^2\Phi$. This corresponds to estimating $\partial_t\mathcal{L}P$.
A similar result can be obtained for the radial derivatives, using the equation satisfied by the profile. But this is not enough to estimate the standard spatial derivatives, which involve variations of angles. Hence we make the following extra assumption on the behaviour of $\mathcal{L}\Phi$ on average,
\begin{equation}
   \label{eq:decay.space}\int_{\mathbb{S}^{N-1}}|\mathcal{L}\Phi(r\theta)-\mathcal{L}\Phi(r\theta -s\varphi)|\,d\theta\le c\delta r^{-N-\beta}\quad \text{for some } \beta>0,
\end{equation}
whenever $0<s<\delta\le1$, $\varphi\in\mathbb{S}^{N-1}$, $r\ge1$,
which roughly speaking means estimating $\nabla \mathcal{L}P$ on average.
The required smoothness for the function $f$ will be given in  terms of a topology adapted to the scaling of the equation, through the H\"older spaces $C^\alpha_\sigma(Q)$ defined in Section~\ref{sec:distance}.
\begin{Theorem}
\label{thm:regularity}  Assume that  the profile $\Phi$  of the fundamental solution for the operator $\partial_t+\mathcal{L}$ satisfies~\eqref{eq:decay.space} for some $\beta>0$.
If $u$ is a very weak solution  to problem~\eqref{eq:main} with $f\in C^\alpha_\sigma(Q)\cap L^\infty(Q)$ for some $0<\alpha<\min\{\beta,\,\sigma,\,1\}$, then $u\in C^\alpha_\sigma(Q)$.
\end{Theorem}

A slight modification of the proof of this result allows to improve the regularity of the solution at each point where  $f$ is more regular, which is stated in Theorem~\ref{th:regularity.improved}. This will be used in a separate work to prove regularity for the nonlinear equation $\partial_t u+\mathcal{L}\phi(u)=0$; see~\cite{VPQR} for the case of the fractional Laplacian.

The key hypothesis~\eqref{eq:decay.space}  holds, with  $\beta=\sigma$, in the  important special case of stable operators given by sums of fractional Laplacians (of order $\sigma$) of smaller dimensions
\begin{equation}
  \label{muchoslaplas}
  \mathcal{L}=\sum_{k=1}^M(-\Delta_{x^k})^{\sigma/2},\qquad x^k\in\mathbb{R}^{n_k},\quad \sum_{k=1}^Mn_k=N,\quad 1\le M\le N,
\end{equation}
the simplest example being~\eqref{eq:sumalaplas}; see Section~\ref{sec-suma}.  In this special situation we also improve estimate~\eqref{eq:decay.time} to include the critical value $\beta=\sigma$.

\begin{Theorem}
\label{th:sumalaplas}  Conditions~\eqref{eq:decay.time} and~\eqref{eq:decay.space} with $\beta=\sigma$ hold for operators $\mathcal{L}$ of the form~\eqref{muchoslaplas}.
\end{Theorem}
The proof follows by observing that in this case the derivatives of the kernel $P$ can be estimated by $P$ itself, combined with estimate~\eqref{profile-pruitt}. We conjecture that this property is true for the profile of any stable process.  In the case of isotropic processes, even depending on time, such estimates for the spatial derivatives have been obtained in~\cite{Kulczycki-Ryznar}.
\begin{Corollary}
  \label{cor-sumalaplas}
  Let $\mathcal{L}$ be given by~\eqref{muchoslaplas}. If $u$ is a very weak solution  to problem~\eqref{eq:main} with $f\in C^\alpha_\sigma(Q)\cap L^\infty(Q)$ for some $0<\alpha<\min\{\sigma,\,1\}$, then $u\in C^\alpha_\sigma(Q)$.
\end{Corollary}

The obtention of estimates for the derivatives of heat kernels of stable L\'evy processes has been the subject of intensive research in the last years. To this aim an  auxiliary smoothness scale of Haussdorff-type, which we describe next, was introduced in~\cite{Bogdan-Sztonyk-2007}; see also~\cite{Watanabe-2007}.
A measure $\mu$ on $\mathbb{S}^{N-1}$ is said to be a $\gamma$--measure if there is a constant $c>0$ such that
$$
\mu(B(\theta,r)\cap\mathbb{S}^{N-1})\le cr^{\gamma-1} \quad\text{for all }\theta\in\mathbb{S}^{N-1} \text{ and } 0<r<1/2.
$$
It is easy to see that necessarily $\gamma\le N$ and that any finite measure is at least a 1--measure.
The case $\gamma=N$ holds if and only if $\mu$  is absolutely continuous with respect to the Lebesgue measure and has a density function which is bounded. This does not mean that the measure is comparable to that of the isotropic case, since it may degenerate in some directions.  If $\gamma>1$, the measure has no atoms. If $\gamma<N$, it is singular. For instance, a spectral measure satisfying $d\mu(\theta)=a(\theta)d\theta$ where $a(\theta)$ has a singularity of the form $a(\theta)\sim|\theta-\theta_0|^{\gamma-N}$, $\gamma\in(1,N)$, is a $\gamma$--measure.

When the spectral measure $\mu$ is a $\gamma$--measure with $  \sigma+\gamma-N>0$,
it was proved in~\cite[Lemma 2.7]{Bogdan-Sztonik-Knopova} that
\begin{equation}
\label{eq:estimate.derivative}
|\nabla \mathcal{L}\Phi(x)|\le C(1+|x|^2)^{-(\sigma+\gamma)/2},
\end{equation}
which implies  the estimate on average~\eqref{eq:decay.space} for  $\beta=\sigma+\gamma-N$. Thus, when $\mu$ is an $N$--measure, the only restriction in Theorem~\ref{thm:regularity} is $\alpha<\min\{\sigma,1\}$, as in the case in which $\mathcal{L}=(-\Delta)^{\sigma/2}$ studied in~\cite{VPQR}.

Let us remark that, though the pointwise estimate~\eqref{eq:estimate.derivative} is optimal, the integral version that is derived from it seems far from being so if $\gamma<N$, since it does not take into account that the measure of the set of directions in which the derivatives decay slowly is small. Thus, we get a restriction on $\beta$, which cannot be arbitrarily close to $\sigma$. We believe that this restriction is technical, since it does not appear  in the \lq\lq worst'' case~\eqref{eq:sumalaplas}, for which $\gamma=1$;  see Corollary~\ref{cor-sumalaplas}.

When the right-hand side is standard, $\partial_t u+\mathcal{L}u=g$, with $g\in L^\infty(Q)\cap C^\alpha_\sigma(Q)$, $\alpha<\min\{\sigma,1\}$,  very weak solutions satisfy $\partial_t u,\mathcal{L}u\in C^\alpha_\sigma(Q)$. This was proved through a blowup argument combined with a Liouville type theorem in~\cite{FernandezReal-RosOton-2017}. This result follows from ours whenever hypotheses~\eqref{eq:decay.space} holds for every $\beta <\sigma$, which is the case when the operator is given by~\eqref{muchoslaplas}, or when it comes from an $N$--measure. Let us emphasize that the result in~\cite{FernandezReal-RosOton-2017} holds for general stable operators, without any restriction on the spectral measure. However,  the blowup argument used there requires some regularity of the right-hand side term, which is not available for problem~\eqref{eq:main}. This is in fact the main difficulty in the present work.

It is also worth mentioning the papers~\cite{Felsinger-Kassmann-2013, Kassmann-Schwab-2014}, where the authors show H\"older regularity  when $g=0$ for a class of operators $\mathcal{L}$, which are not necessarily translation invariant, that include the special case~\eqref{eq:sumalaplas}. Though their proof may  perhaps be adapted to consider $g\in L^\infty(Q)$, it assumes that the solution lies in the energy space, and hence cannot be used to deal with solutions of problem~\eqref{eq:main} when $f$ is not smooth enough.

Observe finally that if $f$ depends only on $x$ or only on $t$ then problem~\eqref{eq:main} becomes trivial. In the application to nonlinear problems that we have in mind the right-hand sides that arise depend tipically on both variables.

\noindent\textsc{Organization of the paper. } We start with the discussion of the required estimates for the kernel and its derivatives in Section~\ref{sect:properties.kernel}, devoting a separate subsection to the special case of operators of the form~\eqref{muchoslaplas}. Section~\ref{sec-homogeneous} deals with the homogeneous case, $f=0$, which yields uniqueness also when the right-hand side
 is nontrivial. Finally, we consider the problem with a singular forcing in Section~\ref{sect-complete-problem}, proving existence and regularity.


\section{Properties of the heat kernel}\setcounter{equation}{0}
\label{sect:properties.kernel}

The aim of this section is to obtain estimates for the heat kernel and its derivatives allowing to apply Theorem~\ref{thm:regularity} to some families of stable operators. We start by describing estimates which are valid for general stable operators, and pass then to consider the case of operators of the form~\eqref{muchoslaplas}, for which much better estimates are available.

\subsection{General stable operators}\label{sec-general-kernels}

Taking  Fourier transform in \eqref{eq:operador en polares} we get that the multiplier $m$ of the operator $\mathcal{L}$, defined by $\widehat{\mathcal{L}u}(\xi)=m(\xi)\widehat{u}(\xi)$, satisfies
\begin{equation}
\label{multiplier}
\begin{array}{l}
\displaystyle m(\xi)=
\int_{\mathbb{S}^{N-1}}\int_0^\infty
(1-\cos(\xi\cdot (r\theta)))
\frac{dr}{r^{1+\sigma}}\,d\mu(\theta)=|\xi|^\sigma g(\xi/|\xi|),\quad\text{where}
\\[4mm]
\displaystyle g(\zeta)=c_{N,\sigma}\int_{\mathbb {S}^{N-1}}|\zeta\cdot \theta|^\sigma\,d\mu(\theta), \quad\text{and} \\[4mm]
\displaystyle c_{N,\sigma}=\int_0^\infty
(1-\cos t)
\frac{dt}{t^{1+\sigma}}=\frac{\pi^{1/2}\Gamma(1-\sigma/2)}{2^\sigma\sigma\Gamma((1+\sigma)/2)}.
\end{array}
\end{equation}
In particular $m$ is homogeneous of order $\sigma$ and $m(\xi)\sim|\xi|^\sigma$, since by the finiteness of the measure and the non-degeneracy condition~\eqref{eq.ellipcity} we have
\begin{equation}\label{g-bdd}
\lambda c_{N,\sigma}\le g(\zeta)\le\Lambda c_{N,\sigma}.
\end{equation}
The homogeneity of the multiplier implies the homogeneity of the operator,
$$
\phi(x)=u(\lambda x)\quad \Rightarrow \quad
\mathcal{L}\phi(x)=\lambda^\sigma\mathcal{L}u(\lambda x).
$$

On the other hand, by \cite[Theorem 2.4.3]{Samorodnitsky}, any symmetric stable process $X=\{X_t\}_{t\ge0}$ defined on a probability space $(\Omega, \mathbb{F}, \mathds{P})$  has a characteristic function
$$
E[e^{i\xi\cdot X}]=e^{-m(\xi)},
$$
where $m(\xi)$ is given by the L\'evy-Khintchine formula \eqref{multiplier}. We have therefore a one-to-one correspondence between our family of operators $\mathcal{L}$ and the family of symmetric stable processes $X$. If we now   consider the family of probability measures $\{\mathfrak{m}_t\}_{t\ge0}$ on $\mathbb{R}^N$, such that for every Borel set $A\subset\mathbb{R}^N$
\begin{equation*}
 \label{process-measure}
 \int_A d\mathfrak{m}_t\equiv \mathds{P}\big(\{\omega\in \Omega : X_t(\omega)\in A\}\big),
\end{equation*}
we have that $d\mathfrak{m}_t=P(\cdot,t)dx$, and $P$ satisfies the problem
$$
    \partial_tP+\mathcal{L}P=0\quad \text{in } Q, \qquad
    P(\cdot,0)=\delta\quad\text{in }\mathbb{R}^N.
$$
The density function $P$ is usually known as the \emph{transition probability density}, the \emph{Gauss kernel} associated to $\mathcal{L}$, or the \emph{fundamental solution} for the operator $\partial_t+\mathcal{L}$. The homogeneity of the multiplier $m$ gives that this kernel is self-similar,
\begin{equation}\label{perfildeP}
P(x,t)=t^{-N/\sigma}\Phi(xt^{-1/\sigma}),\qquad \widehat\Phi(\xi)=e^{-m(\xi)}.
\end{equation}
Clearly, since  $m(\xi)\ge c|\xi|^\sigma$ we have $\Phi\in C^\infty(\mathbb{R}^N)$, $0\le \Phi\le \int_{\mathbb{R}^N}e^{-m(\xi)}\le c$, and $\int_{\mathbb{R}^N}\Phi=1$. Moreover it is also easy to see that $\Phi$ is strictly positive.

As mentioned in the Introduction, in the isotropic  case $\mathcal{L}=(-\Delta)^{\sigma/2}$ the profile $\Phi$ of the kernel is radial, with a decay $\Phi(r)\sim r^{-N-\sigma}$ for $r$ large.
In the anisotropic case an estimate like the previous one is not true in general. However, as proved in~\cite{Pruitt-Taylor}, this rate of decay holds on average;  see~\eqref{profile-pruitt}. Following the proof of that paper it is not difficult to obtain a decay estimate, on average, of the derivatives of $\Phi$.
\begin{Theorem}
\label{thm:estimates.derivatives}
 For any $k,\,l\in\mathbb{N}\cup\{0\}$ we have
\begin{equation}\label{eq-est-derivatives}
 \int_{\mathbb{S}^{N-1}}\partial_r^l\mathcal{L}^k\Phi(r\theta)\,d\theta=O(r^{-N-(1+(k-1)_+)\sigma-l})\quad\text{as }r\to\infty.
\end{equation}
\end{Theorem}
\noindent{\it Proof. } Consider first the case $l=0$. Since $\widehat{\mathcal{L}^k\Phi}(\xi)=m^k(\xi)e^{-m(\xi)}$, we may write
$$
\begin{array}{rl}
\displaystyle  \int_{\mathbb{S}^{N-1}}\mathcal{L}^k\Phi(r\theta)\,d\theta&\displaystyle=
\int_{\mathbb{S}^{N-1}}\int_{\mathbb{R}^N}m^k(\xi)e^{-m(\xi)}e^{ir\theta\cdot\xi}\, d\xi d\theta \\ [3mm] &\displaystyle=
\int_{\mathbb{R}^N}m^k(\xi)e^{-m(\xi)}\int_{\mathbb{S}^{N-1}}e^{ir\theta\cdot\xi}\, d\theta d\xi .
\end{array}  $$
  The inner integral is computed, using spherical coordinates, in \cite{Pruitt-Taylor},
$$
\int_{\mathbb{S}^{N-1}}e^{ir\theta\cdot\xi}\, d\theta=c(r
|\xi|)^{\frac{2-N}2}J_{\frac{N-2}2}(r|\xi|),
$$
where $J_{\omega}$ is the Bessel function of the first kind of order $\omega$.
We thus get
$$
\begin{array}{rl}
\displaystyle
\int_{\mathbb{S}^{N-1}}\mathcal{L}^k\Phi(r\theta)\,d\theta&\displaystyle=cr^{\frac{2-N}2}
\int_{\mathbb{R}^N}m^k(\xi)e^{-m(\xi)}|\xi|^{\frac{2-N}2}J_{\frac{N-2}2}(r|\xi|)\, d\xi \\ [4mm]
&\displaystyle=cr^{\frac{2-N}2}
\int_{\mathbb{S}^{N-1}}g^k(\eta)\int_0^\infty e^{-g(\eta)s^\sigma}s^{\frac{N}2+k\sigma}J_{\frac{N-2}2}(rs)\, dsd\eta.
\end{array}
$$
We conclude, using \cite[Lemma 1]{Pruitt-Taylor} and \eqref{g-bdd}, the behaviour
$$
  \int_{\mathbb{S}^{N-1}}\mathcal{L}^k\Phi(r\theta)\,d\theta=\begin{cases}
    O(r^{-N-\sigma})&\text{if } k=0,\\ O(r^{-N-k\sigma})&\text{if } k\ge1.
  \end{cases}
$$
To estimate the usual derivatives we use the equation for the profile,
$$
 \sigma \mathcal{L}\Phi=N\Phi+r\partial_r\Phi,
$$
which differentiated gives, for each $k\ge0$, $l\ge0$,
$$
   \sigma \partial_r^l\mathcal{L}^{k+1}\Phi=(N+k\sigma+l)\partial_r^l\mathcal{L}^k\Phi+r\partial_r^{l+1}\mathcal{L}^k\Phi,
  $$and obtain \eqref{eq-est-derivatives} by induction in $l$.
 \qed

Unfortunately the estimates needed in our regularity arguments throughout this paper require taking absolute value before taking the average.
For the fractional derivatives, the  pointwise estimate
\begin{equation}\label{profile-glowacki}
|\mathcal{L}^k\Phi(r\theta)|\le \Omega_{\beta,k}(\theta)r^{-N-\beta}\quad \text{a.~e. }\theta\in\mathbb{S}^{N-1},
\end{equation}
was obtained in~\cite{Glowacki-1993} for every $r\ge1$, $k\ge1$ and every $0<\beta<\sigma$, where $\Omega_{\beta,k}\in L^1(\mathbb{S}^{N-1})$. This implies in particular a the decay that is  enough for our purposes.

As we have commented upon in the Introduction,  pointwise estimates for $\nabla\mathcal{L}\Phi$ are not available, except for $\gamma$--measures with $\gamma>N-\sigma$, for which we have~\eqref{eq:estimate.derivative}.

\subsection{The sum of  fractional Laplacians in lower dimensions}\label{sec-suma}

We now turn our attention to the interesting model of stable operators \eqref{eq:operador en polares} of the form~\eqref{muchoslaplas}. Our aim is to show that condition~\eqref{eq:decay.space} holds, so that Theorem~\ref{thm:regularity} can be applied.

For the reader's convenience we perform the calculations in detail. We thus consider sums of fractional Laplacians  $(-\Delta_{x^k})^{\sigma/2}$, whose action on functions of $x=(x',x^k,x'')$, $x'\in\mathbb{R}^{n_1}\times\cdots\times\mathbb{R}^{n_{k-1}}$, $x''\in\mathbb{R}^{n_{k+1}}\times\cdots\times\mathbb{R}^{n_{M}}$, is defined   by
$$\begin{array}{l}
\displaystyle(-\Delta_{x^k})^{\sigma/2}u(x)= \\ [4mm]
\qquad\displaystyle c_{n_k,\sigma}\int_{\mathbb{R}^{n_k}}
\left(u(x',x^k,x'')-\frac{u(x',x^k + y^k,x'') + u(x',x^k-y^k,x'')}2 \right)|y^k|^{-n_k-\sigma}\,dy^k.
\end{array}
$$
The normalization constant $c_{n_k,\sigma}$ is chosen so that the symbol of that operator is $
m_k(\xi)=|\xi^k|^\sigma$, see \eqref{multiplier}, and thus the symbol of $\mathcal{L}$ is
$$
m(\xi)=\sum_{k=1}^M|\xi^k|^\sigma, \qquad \xi=(\xi^1,\cdots, \xi^M),\quad \xi^k\in\mathbb{R}^{n_k}.
$$
The spectral measure of $\mathcal{L}$ is
$$
d\mu(\theta)=\sum_{k=1}^M c_{n_k,\sigma}\delta_{\Omega_k}(\theta),\qquad \Omega_k=\mathbb{S}^{N-1}_+\cap \mathbb{R}^{n_k}.
$$
The most relevant case is when $\mathcal{L}$ is the sum of fractional Laplacians of dimension one, cf.~\eqref{eq:sumalaplas},
for which the spectral measure  is
$d\mu(\theta)=c_{1,\sigma}\sum_{j=1}^N\delta_{e_j}(\theta)$, where
$\{e_j\}_{j=1}^N$ is the canonical basis in $\mathbb{R}^{N}$. Actually we have
$$\begin{array}{l}
\displaystyle\int_{\mathbb{S}^{N-1}}\int_0^\infty
\left(u( x)-\frac{u(x + r\theta) + u(x-r\theta)}2 \right)
\frac{dr}{r^{1+\sigma}}\,d\mu(\theta)\\ [4mm]
\qquad =\displaystyle c_{1,\sigma}\sum_{j=1}^N\int_0^\infty
\left(u( x)-\frac{u(x + re_j) + u(x-re_j)}2 \right)
\frac{dr}{r^{1+\sigma}}= \sum_{j=1}^N(-\partial_{x_jx_j}^2)^{\sigma/2}u(x).
\end{array}
$$

The  operator  \eqref{muchoslaplas}  is the infinitesimal generator of the L\'evy process in $\mathbb{R}^N$ given by $X=\{X_t\}_{t\ge0}$, with $X_t = (X^1_t ,\dots, X^M_t )$,  and $X^k_t$ being independent symmetric stable processes in dimension $n_k$.
The kernel associated to these processes has a profile in separated variables,
\begin{equation}\label{eq.sum-lap-ker}
\Phi(x)=\prod_{k=1}^M\Psi_{k}(x^k),
\end{equation}
where $\Psi_{k}$ is the profile of the kernel corresponding to $(-\Delta_{x^k})^{\sigma/2}$.
This kernel is explicit only when $\sigma=1$, $\Psi_k(w)=d_{n_k}(1+|w|^2)^{-\frac{n_k+1}2}$.
In this particular case,  if we let $x$ tend to infinity along one of the axes $x_j\in x^k$ (see notation below), then $\Phi(x)\sim\Psi_k(x^k)\sim|x^k|^{-n_k-1}$. Thus, the first estimate in~\eqref{eq:first-decay} is not satisfied. The same happens for any $\sigma\in(0,2)$. This example motivates the use of estimates on average on $\mathbb{S}^{N-1}$; see~\cite{Pruitt-Taylor}.

The proof of~\eqref{eq:decay.space} when $\mathcal{L}$ is given by~\eqref{muchoslaplas}  relies on an explicit calculation and an estimate of the kernels $\Psi_k$.
\begin{Proposition}
\label{prop:estimate.gradient}
  The profile $\Phi$ in \eqref{eq.sum-lap-ker} satisfies
  \begin{equation}\label{eq:decay-Phis-radial}
    |\mathcal{L}\Phi(x)|,\,|\nabla\Phi(x)|,\,|\nabla\mathcal{L}\Phi(x)|\le c\Phi(x).
  \end{equation}
\end{Proposition}

\noindent{\it Proof. } First of all we observe that each $\Psi_k$ is radial, so that by Fourier transform as in \cite{Pruitt-Taylor}, see also \cite{VPQR}, we have for every $w\in\mathbb{R}^{n_k}$,
\begin{equation}\label{eq:estimates-k}
\begin{array}{l}
\displaystyle c_1(1+|w|^2)^{-\frac{n_k+\sigma}2}\le\Psi_k(w)\le c_2(1+|w|^2)^{-\frac{n_k+\sigma}2},\\ [4mm]
\displaystyle|L_k\Psi_k(w)|\le c(1+|w|^2)^{-\frac{n_k+\sigma}2},\\ [4mm]\displaystyle|\nabla\Psi_k(w)|,\;|\nabla L_k\Psi_k(w)|\le c(1+|w|^2)^{-\frac{n_k+\sigma+1}2}.
\end{array}
\end{equation}
We have denoted $L_k=(-\Delta_{x^k})^{\sigma/2}$. Therefore \eqref{eq:decay-Phis-radial} holds for each factor in the product~\eqref{eq.sum-lap-ker}. We also use the following convention
$$
\begin{array}{c}
\displaystyle x=(x_1,x_2,\cdots,x_N)=(x^1,x^2,\cdots,x^M),\qquad x^j=(x^j_1,x^j_2,\cdots,x^j_{n_j}), \\ [3mm]\displaystyle x_m\in x^j\;\Leftrightarrow\; x_m=x^j_\ell \text{ for some }\; \ell\in\{1,2,\cdots,n_j\}.
\end{array}
$$
We now calculate
$$
\begin{array}{l}
\displaystyle\partial_{x_m}\Phi(x)=\partial_{x_m}\Psi_j(x^j)\prod_{i\ne j}\Psi_i(x^i)=
\frac{\partial_{x_m}\Psi_j(x^j)}{\Psi_j(x^j)}\,\Phi(x)\quad\text{if } x_m\in x^j,\\[4mm]
\displaystyle\mathcal{L}\Phi(x)=\sum_{k=1}^ML_k\Psi_k(x^k)\prod_{i\ne k}\Psi_i(x^i)=
\sum_{k=1}^M\frac{L_k\Psi_k(x^k)}{\Psi_k(x^k)}\,\Phi(x),
\\[4mm]
\displaystyle\partial_{x_m}\mathcal{L}\Phi(x)\displaystyle=\partial_{x_m}\left(\frac{L_j\Psi_j(x^j)}{\Psi_j(x^j)}\right)\Phi(x)+ \sum_{k=1}^M\frac{L_k\Psi_k(x^k)}{\Psi_k(x^k)}\partial_{x_m}\Phi(x) \\ [4mm]
\displaystyle\qquad\quad=\left(\frac{\partial_{x_m} L_j\Psi_j(x^j)}{\Psi_j(x^j)}-L_j\Psi_j(x^j)\frac{\partial_{x_m}\Psi_j(x_j)}{\Psi^2_j(x^j)}+
\frac{\partial_{x_m}\Psi_j(x^j)}{\Psi_j(x^j)}\sum_{k=1}^M\frac{L_k\Psi_k(x^k)}{\Psi_k(x^k)}\right)\,\Phi(x).
\end{array}
$$
Therefore,
$$
\nabla\mathcal{L}\Phi
=\left(\left.\left\{\frac{\nabla L_j\Psi_j}{\Psi_j}-\frac{L_j\Psi_j\nabla\Psi_j}{\Psi^2_j}+\frac{\nabla\Psi_j}{\Psi_j}
\sum_{j=1}^M\frac{L_k\Psi_k}{\Psi_k}\right\}\right|_{j=1}^M\right)\,\Phi.
$$

We conclude that each coefficient of $\Phi$ in the above derivatives is bounded.
\qed

\medskip

\noindent\emph{Remark. } Actually, estimate~\eqref{eq:estimates-k} implies a sharper estimate for the gradient of $\Phi$,
$$
|\nabla\Phi(x)|^2=\sum_{j=1}^M\frac{|\nabla\Psi_j(x^j)|^2}{\Psi_j^2(x^j)}\Phi^2(x)\le
c\sum_{j=1}^M\frac1{(1+|x^j|^2)}\Phi^2(x).
$$
This gives $|\nabla\Phi(x)|\le c|x|^{-N-\sigma-1}$, as in the radial case, provided $|x|$ is large with $|x^j|\sim|x^k|$ for every $j,k$. In the same way $|\nabla\mathcal{L}\Phi(x)|\le c|x|^{-N-\sigma-1}$ for those directions.

\medskip

We have now the ingredients to prove Theorem~\ref{th:sumalaplas}.

\noindent{\it Proof of Theorem~\ref{th:sumalaplas}. }
We estimate the difference within the integral~\eqref{eq:decay.space} by the Mean Value Theorem. Thanks to Proposition~\ref{prop:estimate.gradient} this amounts to estimate $\int_{\mathbb{S}^{N-1}}\Phi( r\theta-\lambda s\varphi)\,d\theta$, $0\le\lambda\le1$, where $\lambda$ may depend on $\varphi$. In order to use now the estimate on average~\eqref{profile-pruitt}, which would conclude the proof, we must check that we can replace $r\theta-\lambda s\varphi$ in the integral by $r\theta$. If this were the case
$$
\begin{array}{rl}
\displaystyle\int_{\mathbb{S}^{N-1}}|\mathcal{L}\Phi(r\theta)-\mathcal{L}\Phi(r\theta -s\varphi)|\,d\theta &\displaystyle=s\int_0^1 \int_{\mathbb{S}^{N-1}}|\nabla\mathcal{L}\Phi(r\theta-\lambda s\varphi)|\,d\theta d\lambda \\ [4mm]
&\displaystyle\le c\delta\int_0^1 \int_{\mathbb{S}^{N-1}}|\Phi(r\theta-\lambda s\varphi)|\,d\theta d\lambda\\ [4mm]
&\displaystyle\le c\delta\int_{\mathbb{S}^{N-1}}|\Phi(r\theta)|\,d\theta d\le c\delta r^{-N-\sigma}.
\end{array}
$$
So it is enough to prove that $\Phi(z)\le c\Phi(x)$ whenever $|x-z|<\delta\le1$ and $|x|>2\delta$. In fact, since $|x-z|^2=\sum_{k=1}^M|x^k-z^k|^2$, we have $|x^k-z^k|<\delta$ for every $k=1,\cdots,M$. If $|x^k|>2\delta$ this implies $|z^k|>\frac{|x^k|}2$, and thus $\Psi_k(z^k)\le \frac{2^{n_k+\sigma}c_2}{c_1}\Psi(x^k)$ by~\eqref{eq:estimates-k}. On the other hand, if $|x^k|\le2\delta\le2$ we have, again by~\eqref{eq:estimates-k}, $\Psi_k(z^k)\le \frac{c_25^{\frac{n_k+\sigma}2}}{c_1}\Psi(x^k)$. The claim is proved by multiplying all the factors in $k$, and so is the theorem.
\qed

\medskip

\noindent\emph{Remark. } For general operators of the form \eqref{eq:operador en polares}, even if we had an estimate like $|\nabla\mathcal{L}\Phi(x)|\le c\Phi(x)$, this would not imply directly~\eqref{eq:decay.space} as in the previous case, where the special form~\eqref{eq.sum-lap-ker} was used, and further investigation would be needed.
\section{The homogeneous problem}\label{sec-homogeneous}

\label{sect-kernel} \setcounter{equation}{0}

We consider in this section problem \eqref{eq:main} with $f\equiv0$ and prove existence and uniqueness of a very weak solution for every initial datum $u_0\in L^p(\mathbb{R}^N)$ for some $1\le p\le \infty$. To define such concept of solution we consider the weighted space $L_\rho=L^1(\mathbb{R}^N,
\rho\,dx)$ with  weight $\rho(x)=(1+|x|)^{-(N+\sigma)}$. We say that $u\in L^1_{loc}((0,\infty): L_\rho)$ is a  {\sl very weak solution} to problem \eqref{eq:main}  with $f\equiv0$ if
\begin{equation}\label{veryweak0}
\int_{Q} u(\partial_t\zeta - \mathcal{L}\zeta)\,dxdt+\int_{\mathbb{R}^N}u_0(x)\zeta(x,0)\,dx=0
\end{equation}
for all $\zeta\in C_0^\infty(\overline Q)$. The introduction of the weighted space $L_\rho$ allows for the term $\int_Q u\mathcal{L}\zeta$ to be well defined, due to the decay of $\mathcal{L}\zeta$. In fact by a classical result on Fourier Analysis, $m(\xi)\sim|\xi|^\sigma$ implies $\mathcal{L}\zeta(x)=O(|x|^{-(N+\sigma)})$ for large $|x|$.

We will show that test functions which are not compactly supported are also admissible,  provided they have a minimal decay at infinity. In order to prove this assertion we will use the formula contained in the next proposition, which follows easily from a direct computation.
\begin{Proposition}
	\label{pro:L-product} For every pair $v,\,w\in L_\rho$,
	$$
	\mathcal{L}(vw)=v\mathcal{L}(w)+w\mathcal{L}(v)-E(v,w),
	$$
	where
	\begin{equation}\label{resto}
	\begin{array}{rl}
	E(v,w)(x)&\displaystyle=\frac12\int_{\mathbb{R}^{N}}
	\Big(v(x+y)-v(x) \Big)\Big(w(x+y)-w(x) \Big)\,d\nu(y) \\ [4mm]
	&\displaystyle+\frac12\int_{\mathbb{R}^{N}}
	\Big(v(x)-v(x-y) \Big)\Big(w(x)-w(x-y) \Big)\,d\nu(y).
	\end{array}
	\end{equation}
\end{Proposition}
Observe that if the measure $\nu$ were symmetric this expression would simplify to
$$
E(v,w)(x)=\int_{\mathbb{R}^{N}}
\Big(v(x+y)-v(x) \Big)\Big(w(x+y)-w(x) \Big)\,d\nu(y),
$$
formula that appears in~\cite{Barrios-Peral-Soria-Valdinoci} for the case of the fractional Laplacian.

\begin{Proposition}\label{pro:moretests}
Let $u$ be a   very weak solution to problem \eqref{eq:main} with $f\equiv0$. Let also $\varphi\in C^\infty(Q)$ be a function that vanishes for $t>t_0$ for some $t_0>0$, and satisfies $\varphi,\,|\nabla\varphi|\le c\rho$ in $\mathbb{R}^N\times[0,t_0]$. Then identity \eqref{veryweak0} holds with $\zeta$ replaced by $\varphi$.
\end{Proposition}

\noindent{\it Proof. }
We  multiply $\varphi$ by a sequence of cut-off functions, use  identity \eqref{veryweak0} with these admissible test functions and  pass to the limit.

Let then $\phi\in C_0^\infty(\mathbb{R}_+)$ be a nonincreasing function such that $\phi\equiv1$ for $0\le s\le1/2$ and $\phi\equiv0$ for $s\ge1$,  and define the function $\phi_R(x)=\phi(R^{-1}|x|)$.
We  are done if we show that
\begin{equation}\label{new-test}
\lim_{R\to\infty}\int_{\mathbb{R}^N}u\mathcal{L}(\varphi\phi_R)=\int_{\mathbb{R}^N}u\mathcal{L}\varphi
\end{equation}
for each fixed time $0<t<t_0$. In order to do that we need to compute the action of $\mathcal{L}$ on the product $\varphi\phi_R$.
 Since the bilinear form $E(\varphi,\phi_R)$ only involves products of differences, see~\eqref{resto}, using the same proof as in~\cite{Barrios-Peral-Soria-Valdinoci} we obtain
$$
\lim_{R\to\infty}\int_{\mathbb{R}^N}uE(\varphi,\phi_R)=0.
$$
The main point is the hypothesis $\varphi,\,|\nabla\varphi|\le c\rho$.
Recall finally that we have $\|\mathcal{L}\phi\|_\infty\le c(\|\phi\|_\infty+\|D^2\phi\|_\infty)$, so that by homogeneity, and the fact that $u\varphi\in L^1(\mathbb{R}^N)$,
$$
\lim_{R\to\infty}\int_{\mathbb{R}^N}u\varphi\mathcal{L}\phi_R=
\lim_{R\to\infty}R^{-\sigma}\int_{\mathbb{R}^N}u\varphi\mathcal{L}\phi=0.
$$
We therefore get \eqref{new-test}.
\qed

\begin{Theorem}
\label{th:f=0}   If $u_0\in L^p(\mathbb{R}^N)$ for some $1\le p\le \infty$, then  problem \eqref{eq:main} with $f\equiv0$ has a unique
very weak solution. The solution is bounded and $C^\infty$ smooth for every $t>0$ and satisfies the equation in the classical sense.
\end{Theorem}

\noindent{\it Proof. }
Existence follows easily by convolution with the heat kernel, $u=u_0*P$. Thus we deduce the same standard smoothing effect as for the solutions of the fractional heat equation, or even the local heat equation: $u(\cdot,t)\in C^\infty(\mathbb{R}^N)\cap L^q(\mathbb{R}^N)$ for every $p\le q\le \infty$ and any $t>0$, with
$$
\|u(\cdot,t)\|_q\le ct^{-\frac N\sigma(\frac1p-\frac1q)}\|u_0\|_p.
$$
In order to prove uniqueness we just consider the case $u_0\equiv0$, then take $R_0,\,t_0>0$ arbitrary and show that
\begin{equation}
  \label{u=0}
  \int_Qu(x,t)F(x,t)\,dxdt=0 \quad\text{for all }F\in C^\infty_0(\{|x|<R_0,\,0<t<t_0\}).
\end{equation}
We use Hilbert's duality method by considering as test function in the definition of very weak solution  the unique solution $\varphi$ to the nonhomogeneous backward problem
\begin{equation}
\label{eq:dual.problem}
\left\{
\begin{array}{ll}
\partial_t\varphi-\mathcal{L}\varphi=F,\quad &x\in\mathbb{R}^N,\,0<t<t_0,
\\ [4mm]
\varphi = 0,\quad &x\in \mathbb{R}^N,\,t\ge t_0.
\end{array}
\right.
\end{equation}
This would yield \eqref{u=0} once we check that $\varphi$ is a good test function.
Though the fact that $\mathcal{L}$ is a non-local operator implies that $\varphi$ does not have compact support in $\overline Q$, Proposition~\ref{pro:moretests} allows to use it as a test function provided $\varphi\le c\rho$ and $|\nabla \varphi|\le c\rho$.

Using Duhamel's formula, a solution to~\eqref{eq:dual.problem} can be written using the heat kernel $P$ in the form
$$
\varphi(x,t)=\int_0^{t_0-t}\int_{\mathbb{R}^N}P(x-y,t_0-t-s)F(y,s)\,dyds.
$$
By Young's inequality we have $\|\varphi(\cdot,t)\|_\infty\le\|F\|_\infty$. On the other hand, since $F$ has compact support, taking $|x|>2R_0$, we have, using the self-similar form of $P$ and \eqref{profile-pruitt},
$$
\begin{array}{rl}
\displaystyle|\varphi(x,t)|&\le\displaystyle\|F\|_\infty\int_0^{t_0-t}\int_{|y|<R_0}P(x-y,t_0-t-s)\,dyds \\ [4mm]
&\le\displaystyle \|F\|_\infty\int_0^{t_0-t}\int_0^{R_0}\int_{\mathbb{S}^{N-1}}P(x-r\theta,s)r^{N-1}\,drd\theta ds \\ [4mm]
&\le\displaystyle c\|F\|_\infty\int_0^{t_0-t}\int_0^{R_0}s^{-\frac N\sigma}\left(|x-r\theta|s^{-\frac1\sigma}\right)^{-N-\sigma}r^{N-1}\,drds \\ [4mm]
&\le\displaystyle c\|F\|_\infty(t_0-t)^2\int_0^{R_0}\frac{r^{N-1}}{|x-r\theta|^{N+\sigma}}\,dr \\ [4mm]
&\le\displaystyle \frac{c\|F\|_\infty(t_0-t)^2R_0^N}{(1+|x|^2)^{\frac{N+\sigma}2}},
\end{array}
$$
since $|x-r\theta|>|x|/2>c(1+|x|^2)^{\frac12}$. In the same way we estimate $|\nabla\varphi|$, this time in terms of $\|\nabla F\|_\infty$. We end the proof as follows: use identity~\eqref{veryweak0} with $u_0\equiv0$ and test function $\zeta=\varphi$ solution to problem~\eqref{eq:dual.problem}, which gives~\eqref{u=0} and thus $u\equiv0$.
\qed

\section{The problem with reaction}
\label{sect-complete-problem} \setcounter{equation}{0}

We consider here the Cauchy problem  \eqref{eq:main} with a nontrivial right-hand side. Since the equation is linear, thanks to the previous section we may assume without loss of generality that $u_0\equiv0$.  We  define a  {\sl very weak solution} to problem \eqref{eq:main} (with $u_0\equiv0$) as a function   $u\in L^1_{loc}([0,\infty): L_\rho)$ such that
\begin{equation}\label{veryweak}
\int_{Q} u\partial_t\zeta =\int_{Q} (u-f)\mathcal{L}\zeta\quad \text{for all } \zeta\in C_0^\infty(\overline Q).
\end{equation}

\subsection{$\sigma$--parabolic distance}\label{sec:distance}

We introduce now a topology adapted to the equation, in terms of which the estimates are easier to write.  This notation has already been used in the literature; see for instance~\cite{VPQR}.
In order to reflect the different influence of the variables in the equation, we use a {\em $\sigma$-parabolic ``distance''} $|Y_1-Y_2|_\sigma$ between points $Y_1,\,Y_2\in Q$, derived from the {\em $\sigma$-parabolic ``norm''} defined by
\begin{equation*}
\label{sigma-distance}
|Y|_\sigma:=\Big(|x|^2+|t|^{2/\sigma}\Big)^{1/2}=|t|^{1/\sigma}(|z|^2+1)^{1/2},\qquad
Y=(x,t)\in\mathbb{R}^N\times\mathbb{R},\quad z=x|t|^{-1/\sigma}.
\end{equation*}
We clearly have $ |Y|\le |Y|_\sigma^{\min\{\sigma,\,1\}}$.
The  H\"older space  $C^\alpha_\sigma(Q)$, $\alpha\in(0,\nu)$, will consist of functions $u$ defined in $Q$ such that for some constant $c>0$
$$
|u(Y_1)-u(Y_2)|\le c|Y_1-Y_2|_\sigma^\alpha \quad \text{ for every }Y_1,\,Y_2\in Q.
$$

The $\sigma$-parabolic ball  is defined as  $B_R:=\{Y\in\mathbb{R}^{N+1}\,:\,|Y|_\sigma<R\}$. It is also useful to write each point $x\in\mathbb{R}^N$ in polar coordinates, $x=r\theta$, $r=|x|$, $\theta=x/r$. In that way, to each point $Y=(x,t)=(r\theta,t)\in\mathbb{R}^N\times\mathbb{R}$ we associate the point $\widetilde Y=(r,t)\in\mathbb{R}^2_{+}$  and write, by abuse of notation,  $Y=(\widetilde Y,\theta)$. Observe that, again abusing notation,
$$
|Y|_\sigma=|\widetilde Y|_\sigma=(r^2+|t|^{2/\sigma})^{1/2}.
$$

Let us also consider the ball $\widetilde B_R:=\{\widetilde Y\in\mathbb{R}^2_+\,:\,|\widetilde Y|_\sigma<R\}$.
We  have  that the integrals in $\sigma$-parabolic balls  can be decomposed as
$$
\begin{array}{rl}
\displaystyle\int_{B_R}w(Y)\,dY&\displaystyle
=\int_{\mathbb{S}^{N-1}}\int_{\{r^2+|t|^{2/\sigma}<R^2\}}w(r,t,\theta)r^{N-1}\,drdtd\theta \\ [4mm]
&\displaystyle=\int_{\mathbb{S}^{N-1}}\int_{\widetilde B_R}w(\widetilde Y,\theta)r^{N-1}\,d\widetilde Yd\theta.
\end{array}
$$

For instance, by using the change of variables
\begin{equation}\label{changevariables}
s=r|t|^{-1/\sigma},\qquad\rho=(r^2+|t|^{2/\sigma})^{1/2},
\end{equation}
we can obtain
$$
\begin{array}{rl}
\displaystyle\int_{B_{R}}g(|Y|_\sigma)\,dY&\displaystyle=2\omega_N\int_0^\infty\int_0^R g(\rho)\left(\frac{s\rho}{(s^2+1)^{1/2}}\right)^{N-1}
\frac{\sigma\rho^\sigma}{(s^2+1)^{\frac{\sigma+1}2}}\,d\rho ds\\ [3mm]
&\displaystyle=c\int_0^Rg(\rho)\rho^{N+\sigma-1}\,d\rho.
\end{array}
$$
In particular, the volume of the ball  $B_R$  is proportional to
$R^{N+\sigma}$.

We finally write, in terms of the new distance, the estimates  for the time derivatives of the Gauss kernel $P$ that can be deduced from the decay estimates for the profile $\Phi$, see Section~\ref{sec-general-kernels}. We use self-similarity and the fact that  $\Phi$ is bounded, so any estimate $cr^{-N-\beta}$ for $r$ large can be written as $c(1+r^2)^{\frac{-N-\beta}2}$ for $r>0$.
For every $k\in\mathbb{N}$ it holds
\begin{equation}\label{general-estimate-kernel1}
\begin{array}{rl}
\displaystyle\int_{\mathbb{S}^{N-1}}|\partial_t^kP(r\theta,t)|\,d\theta
&\displaystyle=\int_{\mathbb{S}^{N-1}}|\mathcal{L}^kP(r\theta,t)|\,d\theta \\ [4mm]
&\displaystyle=t^{-\frac N\sigma-k}\int_{\mathbb{S}^{N-1}}|\mathcal{L}^k\Phi(r\theta t^{-\frac1\sigma})|\,d\theta \\ [4mm]
&\displaystyle\le ct^{-\frac N\sigma-k}(1+r^2t^{-\frac2\sigma})^{\frac{-N-\beta}2} \\ [4mm]
&\displaystyle\le c t^{-k+\frac\beta\sigma}|\widetilde{Y}|_\sigma^{-N-\beta},
\end{array}
\end{equation}
for every $0<\beta<\sigma$, $\widetilde Y=(r,t)$, $r,t>0$. If $k=0$ it is also true with $\beta=\sigma$.

\subsection{A cancelation property}\label{sec-cancel}

We next show a cancelation property for $\mathcal{L}P$ crucial in later regularity arguments.
\begin{Theorem}\label{pro:average0} Let $P$ be the Gauss kernel \eqref{perfildeP}.
Then for every $0<a<b$,
\begin{equation}
  \label{int-A=0}
\int_{B_{a,b}^+}\mathcal{L}P(Y)\,dY=0,
\end{equation}
where $B_{a,b}^+=\{a<|Y|_\sigma<b,\,t>0\}$.
\end{Theorem}

\noindent{\it Proof. } Using as before the change of variables \eqref{changevariables}
we get,
$$
\begin{array}{rl}
\displaystyle\int_{B_{a,b}^+}\mathcal{L}P(Y)\,dY&=
\displaystyle-\int_{B_{a,b}^+}\partial_tP(Y)\,dY\\ [3mm]&\displaystyle= \frac1\sigma\int_{\mathbb{S}^{N-1}}
\int_{\widetilde B_{a,b}^+} t^{-\frac{N+\sigma}\sigma}(N\Phi(r\theta t^{-\frac1\sigma})+rt^{-1/\sigma}\theta\nabla \Phi(r\theta t^{-1/\sigma}))r^{N-1}\,drdtd\theta\\ [3mm]
&\displaystyle=
\log(
b/a)\int_{\mathbb{S}^{N-1}}\int_0^\infty s^{N-1} (N\Phi(s\theta)+s\theta\partial_s \Phi(s\theta))\,ds d\theta \\ [3mm]
&\displaystyle=
\log(
b/a) \int_{\mathbb{S}^{N-1}}\int_0^\infty \partial_s(s^N \Phi(s\theta))\,dsd\theta\\ [3mm]
&\displaystyle=\log(
b/a) \lim_{R\to\infty}R^N \int_{\mathbb{S}^{N-1}}\Phi(R\theta))\,d\theta
=0.
\end{array}
$$
The last limit uses the behaviour \eqref{profile-pruitt}.
 \qed

\subsection{Existence}

We formally write the solution using Duhamel's formula:
$$
u(x,t)=\int_0^t\int_{\mathbb{R}^N}P(x-\overline x,
t-\overline t) \mathcal{L}f(\overline x,\overline t)\,d\overline xd\overline t,\quad (x,t)\in Q.
$$
Now integrate by parts  and consider the integral in principal value sense (in $\sigma$--parabolic topology). We  prove that what we obtain is in fact the unique solution to our problem.
\begin{Theorem}
\label{th:duhamel.very.weak}   If $f\in C^\alpha_\sigma(Q)\cap L^\infty(Q)$ for some $0<\alpha<1$, then  the function
\begin{equation}\label{eq:duhamel1}
u(x,t)=\lim_{\varepsilon\to0}\int_{\Omega_{\varepsilon}(x,t)}\mathcal{L}P(x-\overline x,
t-\overline t) f(\overline x,\overline t)\,d\overline xd\overline t,
\end{equation}
where $\Omega_{\varepsilon}(x,t)=\{\varepsilon^2<|\overline x-x|^2+|\overline t-t|^{2/\sigma},\,0<\overline t<t \}$, is the unique
very weak solution of problem~\eqref{eq:main} with $u_0\equiv0$, which is moreover bounded.
\end{Theorem}

\noindent{\it Proof. } Uniqueness follows from the previous section. Let us show that  the function in \eqref{eq:duhamel1}  is well defined.
Let $Y=(x,t)\in Q$ be fixed and take $0<\varepsilon<t^{1/\sigma}/2$. We decompose the integral  as
$$
\int_{\Omega_{\varepsilon}(x,t)}=\int_{\Omega_{\varepsilon}(x,t)-\Omega_{t^{1/\sigma}}(x,t)}+
\int_{\Omega_{t^{1/\sigma}}(x,t)}=I_1+I_2.
$$
The cancellation property \eqref{int-A=0} implies
$$
\int_{\Omega_{\varepsilon}(x,t)-\Omega_{t^{1/\sigma}}(x,t)}\mathcal{L}P(x-\overline x,
t-\overline t) \,d\overline xd\overline t=0.
$$
Therefore,  the   H\"older regularity  of $f$ together with estimate \eqref{general-estimate-kernel1} with $k=1$ and $\beta<\sigma$ allow us to estimate the inner integral,
$$
\begin{array}{rcl}
\displaystyle|I_1|&=&
\displaystyle
\left|\int_{\Omega_{\varepsilon}(Y)-\Omega_{t^{1/\sigma}}(Y)}\mathcal{L}P(Y-\overline Y)f(\overline Y)\,d\overline Y\right|
\\ [4mm]&=&
\displaystyle
\left|\int_{\Omega_{\varepsilon}(Y)-\Omega_{t^{1/\sigma}}(Y)}\mathcal{L}P(Y-\overline Y)\Big(f(\overline Y)-f(Y)\Big)\,d\overline Y\right|
\\ [4mm]
&\le&\displaystyle [f]_{C^\alpha}\int_{|Y-\overline Y|^\sigma_\sigma<t}|\mathcal{L}P(Y-\overline Y)||\overline Y- Y|_\sigma^{\alpha}\,d\overline Y
\\ [4mm]
&=&\displaystyle c\int_{\mathbb{S}^{N-1}}\int_{|\widetilde Y|^\sigma_\sigma<t}|\mathcal{L}P(\widetilde Y,\theta)||\widetilde Y|_\sigma^{\alpha} r^{N-1}\,d\widetilde Yd\theta \\ [4mm]
&\le&\displaystyle  c\int_{|\widetilde Y|^\sigma_\sigma<t}\tau^{-1+\beta/\sigma}|\widetilde Y|_\sigma^{-N-\beta+\alpha}r^{N-1}\,d\widetilde Y\\ [4mm]
&=&\displaystyle c \int_0^\infty\dfrac{s^{N-1}}{(s^2+1)^{\frac{N+\beta}2}}\,ds\,\int_0^{t^{1/\sigma}}\rho^{\alpha-1}\,d\rho
=ct^{\alpha/\sigma}.
\end{array}
$$
We have put $Y-\overline Y=(\widetilde Y,\theta)=(r,\tau,\theta)$, integrated in the sphere and then used the change of variables \eqref{changevariables}.

We next prove that the outer integral is bounded by using the boundedness of $f$. Here we integrate first in the sphere, then in the radial variable and finally in time.
$$
\begin{aligned}
\displaystyle|I_2|\le& \displaystyle \|f\|_\infty
\int_{\Omega_{t^{1/\sigma}}(x,t)}|\mathcal{L}P(Y-\overline Y)|\,d\overline Y
\\
\le&\displaystyle c
\int_{\mathbb{S}^{N-1}}\int_{|\widetilde Y|^\sigma_\sigma>t}|\mathcal{L}P(\widetilde Y,\theta)|r^{N-1}\mathds{1}_{\{0<\tau<t\}} \,d\widetilde Yd\theta \\
\le&\displaystyle
c\int_0^t\int_{\sqrt{t^{2/\sigma}-\tau^{2/\sigma}}}^\infty\tau^{-1+\beta/\sigma}
r^{-\beta-1}\,drd\tau
\\
=&\displaystyle
c\int_0^t\tau^{-1+\beta/\sigma}(t^{2/\sigma}-\tau^{2/\sigma})^{-\beta/2}\,d\tau
=c.
\end{aligned}
$$

This also gives that the solution is bounded for every bounded interval of times. The fact that $u$ is a very weak solution is immediate.  \qed

\subsection{H\"older regularity}\label{sec-regularity}

We study here the regularity of the function given by formula \eqref{eq:duhamel1}, using the notation $Y=(x,t)\in Q$,
\begin{equation*}
  \label{eq:g}
  u(Y)=\int_{\mathbb{R}^{N+1}}A(Y-\overline Y)\mathds{1}_{\{0<\overline t<t\}}f(\overline Y)\,d\overline Y,
\end{equation*}
where $A=\mathcal{L}P$. We omit the principal value sense of the integral for simplicity.

\begin{figure}[ht]
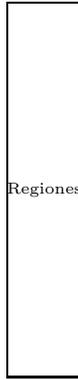

\begin{center}
\pgfuseimage{fig1}\end{center}
\caption{Integration regions for $Y=(r\theta,t)$ for each $\theta\in\mathbb{S}^{N-1}$ fixed.}
\label{fig:integration.regions}
\end{figure}

\noindent{\it Proof of Theorem \ref{thm:regularity}. }  Let $Y_1=(x_1,t_1),\,Y_2=(x_2,t_2)\in Q$ be two points with $|Y_1-
Y_2|_\sigma=h>0$ small, and assume for instance $t_1\ge t_2$. By substracting $f(Y_1)$ to $f$ we may assume without loss of generality that $f(Y_1)=0$. We must estimate the difference
\begin{equation*}\label{eq:g-g0}
\begin{aligned}
|u(Y_1)-u(Y_2)|&\displaystyle=\left|\int_{\mathbb{R}_+^{N+1}}\Big(A( Y_1-
\overline Y)\mathds{1}_{\{\overline t<t_1\}}-A( Y_2- \overline Y)\mathds{1}_{\{\overline t<t_2\}}\Big)f(\overline Y)\,d\overline Y\right| \\
&\displaystyle=\Big|\int_0^{t_1}\int_{\mathbb{R}^{N}}\Big(A(Y)\mathds{1}_{\{t>0\}}-A(Y-Y_3)\mathds{1}_{\{t>t_3\}}\Big)f(Y_1-Y)\,d Y\Big|.
\end{aligned}
\end{equation*}
We have made the change of variables $Y=Y_1-\overline Y=(x,t)$ and put $Y_3=Y_1-Y_2=(x_3,t_3)$, so that $|Y_3|_\sigma=h$, $t_3=t_1-t_2\in[0,h^\sigma]$.  Observe that $|f(Y_1- Y)|\le c|Y|_\sigma^\alpha$.  We decompose $Q_1=\mathbb{R}^N\times(0,t_1)$ into
three regions, depending on the sizes of $| x|$ and $t$, see
Figure~\ref{fig:integration.regions},
$$
\int_{Q_1}=\int_{C_{h}}+\int_{S_h}+\int_{D_h}.
$$

\noindent (i) \emph{The small ``semiball''} $C_h=B^+_{\rho h}=\{|Y|_\sigma<\rho h\}\cap\{t>0\}$, where $\rho>2$ is a constant to be fixed later.   We take $h$ small enough ($(\rho h)^\sigma<\min\{t_1,1\}$) so that
$C_{h}\subset Q_1$.   The difficulty in this region is the
non-integrable singularity of $A( Y)$ at $ Y=0$, which is to be compensated by the regularity of $f$.  We first have, repeating the computations  of the proof of Theorem~\ref{th:duhamel.very.weak},
$$
\int_{C_{h}}\left|A( Y)\right||f(Y_1-Y)|\,d Y
\le\int_{C_{h}}\left|A( Y)\right||Y|_\sigma^\alpha\,d Y
\le
ch^{\alpha}.
$$
As to the second term in $\int_{C_{h}}$, we use the cancelation property \eqref{int-A=0} in order to counteract the singularity at $Y=Y_3$. Thus, taking $\rho\ge\max\{2,2^{1/\sigma}\}$ we have
$$
B^+_{h}(Y_3)=\{|Y-Y_3|_\sigma<h,\,t>t_3\}\subset C_{h}\subset B^+_{\rho^2h}(Y_3),
$$
so that,
\begin{equation}\label{integral-I1}
\begin{aligned}
\left|\int_{C_{h}}A(Y-\right.&\left. Y_3)\mathds{1}_{\{t>t_3\}}f(Y_1-Y)\,d Y\right|
\\
\le&\int_{C_{h}}\left|A(Y-Y_3)\mathds{1}_{\{t>t_3\}}\right|\,
\Big|f(Y_1-Y)-f(Y_1-Y_3)\Big|\,d Y\\
&+\left|f(Y_1-Y_3)\right|\int_{C_{h}-B_{h}(Y_3)}\left|A(Y-Y_3)\mathds{1}_{\{t>t_3\}}\right|\,d Y
\\
\le& c\int_{B^+_{\rho^2 h}(Y_3)}\left|A(Y-Y_3)\right|\,
\left|Y-Y_3\right|^\alpha\,d Y\\
&+\displaystyle
ch^\alpha\int_{B_{\rho^2 h}^+(Y_3)-B_{h}^+(Y_3)}\left|A(Y-Y_3)\right|\,d Y= I_1+I_2.
\end{aligned}
\end{equation}
The first integral satisfies again $I_1\le ch^\alpha$.
We now show that $I_2$ can be controlled since we are far from the singularity. Putting $Z=Y-Y_3$, and using as always the notation in polar coordinates/time, $Z=(s\varphi,\tau)=(\widetilde Z,\varphi)$, we have
$$
\begin{aligned}
I_2&\displaystyle\le c h^{\alpha}\int_{\{\rho h<|Z|_\sigma<\rho^2 h,\,\tau>0\}}|A(Z)|\,dZ \\
&\displaystyle=
c h^{\alpha}\int_{\{\rho h<|\widetilde Z|_\sigma<\rho^2 h,\,\tau>0\}}\int_{\mathbb{S}^{N-1}}|A(\widetilde Z,\varphi)|\,d\varphi d\widetilde Z \\
&\displaystyle\le c h^{\alpha}\int_{\{\rho h\le|\widetilde Z|_\sigma\le\rho^2 h,\,\tau>0\}}\tau^{-1+\frac\beta\sigma}|\widetilde Z|^{-N-\beta}_\sigma\,d\widetilde Z.
\end{aligned}
$$
Changing $(s,\tau)\to(w,\xi)$ following \eqref{changevariables} we end up with the estimate
$$
I_2\le c h^{\alpha}\int_0^\infty\frac{w^{N-1}}{(w^2+1)^{\frac{N+\beta}2}}\,dw\int_{\rho h}^{\rho^2 h}\frac{d\xi}\xi=ch^\alpha.
$$

\noindent (ii) \emph{Outside the ball $B_{\rho h}$ for small
times,
$S_h=\{ |Y|_\sigma\ge\rho h,\;0<t<(\rho h/2)^\sigma\}$. } Since in this region we have $|Y|_\sigma\le \rho_1| Y -Y_3|_\sigma$ for some positive constant $\rho_1$ depending only on $\sigma$, both integrals in $S_h$ are of the same order
$$
\begin{array}{l}
\displaystyle\int_{S_h}\left(|A(Y)|+|A(Y-Y_3)|\mathds{1}_{\{t>t_3\}}\right)\,|f(Y_1-Y)|\,d Y \\ [4mm]
\quad\quad\displaystyle\le
c \int_0^{(\frac{\rho h}2)^\sigma}\int_0^\infty t^{-1+\frac\beta\sigma}(r^2+t^{\frac2\sigma})^{\frac{-N-\beta+\alpha}2}r^{N-1}\,drdt \\ [4mm]
\quad\quad\displaystyle\le c \int_0^{(\frac{\rho h}2)^\sigma} t^{-1+\frac\alpha\sigma}\,dt\int_0^\infty\frac{w^{N-1}}{(w^2+1)^{\frac{N+\beta-\alpha}2}}\,dw
=ch^\alpha.
\end{array}
$$
Notice that   $\alpha<\sigma$ so  the last integral is convergent.

\noindent (iii) \emph{Outside the ball $B_{\rho h}$ for not so small
times, $D_h=\{|Y|_\sigma\ge\rho h,\;(\rho h/2)^\sigma<t<t_1\}$. } Since in that set it is $t>t_3$, we have
$$
A(Y)-A(Y-Y_3)\mathds{1}_{\{t>t_3\}}=A(Y)-A(Y-Y_3),
$$
and there will be some
cancellation. We put $Y_3^*=(0,t_3)$  and decompose this difference as
$$
\displaystyle|A(Y)-A(Y-Y_3)|\le |A(Y)-A(Y-Y_3^*)|+|A(Y-Y_3^*)-A(Y-Y_3)|.
$$
For the first term,
$$
|A(Y)-A(Y-Y_3^*)|\le h^\sigma||\partial_t
A(Z_1)|,
$$
where $Z_1=(x,t-\lambda t_3)$ for some $\lambda\in(0,1)$.
Observe also  that $t-\lambda t_3\sim t$. Using now~\eqref{general-estimate-kernel1} we have, denoting as always, $Y=(r\theta,t)$,
$$
\begin{array}{l}
\displaystyle\int_{D_h}|A(Y)-A(Y-Y_3^*)||f(Y_1-Y)|\,dY \\ [4mm]
\displaystyle\quad\le h^\sigma\int_{D_h}|\partial_t
A(Z_1)||f(Y_1-Y)|\,dY \\ [4mm]
\displaystyle\quad\le ch^\sigma\int_{(\rho h/2)^\sigma}^{t_1}(t-\lambda t_3)^{-\frac N\sigma-2}\int_0^\infty\int_{\mathbb{S}^{N-1}}|\mathcal{L}^2\Phi(r\theta(t-\lambda t_3)^{-\frac1\sigma})|(r^2+t^{\frac2\sigma})^{\frac\alpha2}r^{N-1}\,d\theta drdt \\ [4mm]
\displaystyle \quad\le ch^\sigma\int_{(\rho h/2)^\sigma}^{t_1}t^{-2+\frac\beta\sigma}\int_0^\infty
(r^2+t^{\frac2\sigma})^{\frac{-N-\beta+\alpha}2}r^{N-1}\, drdt\\ [4mm]
\displaystyle \quad= ch^\sigma\int_{(\rho h/2)^\sigma}^{t_1}t^{-2+\frac\alpha\sigma}\int_0^\infty
\frac{w^{N-1}}{(w^2+1)^{\frac{N+\beta-\alpha}2}}\,dwdt\le ch^{\alpha}.
\end{array}
$$
The last integral is convergent  provided $\alpha<\beta$.
We now estimate the spatial difference,
$$
\begin{aligned}
\int_{D_h}|A(Y-Y_3^*)-&A(Y-Y_3)||f(Y_1-Y)|\,dY  \\
&\le c\int_{(\rho h/2)^\sigma}^{t_1}(t- t_3)^{-\frac N\sigma-1}\int_0^\infty I(r,t)(r^2+t^{\frac2\sigma})^{\frac\alpha2}r^{N-1}\, drdt
\end{aligned}
$$
where
$$
I(r,t)=\int_{\mathbb{S}^{N-1}}
|\mathcal{L}\Phi(r\theta(t- t_3)^{-\frac1\sigma})-\mathcal{L}\Phi((r\theta-s\varphi)(t- t_3)^{-\frac1\sigma})|\,d\theta,\qquad s\varphi=x_3.
$$
Using hypothesis~\eqref{eq:decay.space},
$$
I(r,t)\le
cs(t-t_3)^{-\frac1\sigma}(r(t-t_3)^{-\frac1\sigma})^{-N-\beta}\le
cht^{\frac{N+\beta-1}\sigma}(r^2+t^{\frac2\sigma})^{-\frac{N+\beta}2}.
$$
Therefore
$$
\begin{array}{l}
\displaystyle\int_{D_h}|A(Y-Y_3^*)-A(Y-Y_3)||f(Y_1-Y)|\,dY \\ [4mm] \displaystyle\le ch\int_{(\rho h/2)^\sigma}^{t_1}\int_0^\infty t^{-1-\frac{1-\beta}\sigma}
(r^2+t^{\frac2\sigma})^{\frac{-N-\beta+\alpha}2}r^{N-1}\, drdt\\ [4mm]
\displaystyle = ch\int_{(\rho h/2)^\sigma}^{t_1}t^{-1-\frac{1-\alpha}\sigma}\int_0^\infty
\frac{w^{N-1}}{(w^2+1)^{\frac{N+\beta-\alpha}2}}\, dwdt\le ch^{\alpha}.
\end{array}$$
As before we need $\alpha<\beta$. The proof is finished.
\qed

We end with a modification of the previous proof by assuming that the datum $f$ is $C^{\alpha+\epsilon}$ H\"older continuous at some point $Y_1$ and only $C^\alpha$, but with a small coefficient, at the rest of the points, thus getting $C^{\alpha+\epsilon}$ regularity at $Y_1$.
\begin{Theorem}
\label{th:regularity.improved} In the hypotheses of Theorem~\ref{thm:regularity}, assume moreover
that there exist  $c>0$, $\delta_0>0$ and $\epsilon>0$,
$\alpha+\epsilon<\max\{\sigma,1\}$,  such that
\begin{eqnarray}
\label{eq:local.holder.condition2}
|f(Y)-f(Y_1)|\le c|Y-Y_1|_\sigma^{\alpha+\epsilon},\\ [3mm]
\label{eq:local.holder.condition}
|f(Y)-f(\overline Y)|\le c\delta^{\epsilon}\,|Y-\overline Y|_\sigma^{\alpha},
\end{eqnarray}
for all $0<\delta<\delta_0$, $Y,\overline Y\in
B_\delta(Y_1)$. Then,
\begin{equation*}
\label{eq:g.alpha+eps} |u(Y)-u(Y_1)|\le
c'|Y-Y_1|_\sigma^{\alpha+\epsilon},
\end{equation*}
for all $Y\in
B_{\delta_0/2}(Y_1)$, where $c'$ depends on $c$.
\end{Theorem}

\noindent{\it Proof. }  Since $f$ is bounded, condition~\eqref{eq:local.holder.condition2} holds for every $Y\in Q$. This is enough to make all the estimates used to prove Theorem~\ref{thm:regularity}  work, yielding terms which
are $O(h^{\alpha+\epsilon})$, except that for the integral $I_1$ in
\eqref{integral-I1}. To estimate this term, take $\rho h<\delta_0$ and observe that~\eqref{eq:local.holder.condition} gives
$$
\begin{array}{rl}
I_1&=\displaystyle\int_{C_h}|A(Y-Y_3)\mathds{1}_{\{t>t_3\}}
\Big|f(Y_1-Y)-f(Y_1-Y_3)\Big|\,dY \\ [4mm]
&\displaystyle\le
ch^{\epsilon}\int_{C_h}\frac1{|Y-Y_3|_\sigma^{N+\sigma}}|Y-Y_3|_\sigma^\alpha\,d Y\le
ch^{\alpha+\epsilon}.
\end{array}
$$
\qed

This theorem will be used somewhere else to study the regularity of solutions to nonlinear nonlocal equations.


\

\noindent{\large \textbf{Acknowledgments}}

\noindent All authors supported by projects MTM2014-53037-P and MTM2017-87596-P (Spain).



\

\noindent\textbf{Addresses:}

\noindent\textsc{A. de Pablo: } Departamento de Matem\'{a}ticas, Universidad
Carlos III de Madrid, 28911 Legan\'{e}s, Spain. (e-mail: arturo.depablo@uc3m.es).

\noindent\textsc{F. Quir\'{o}s: } Departamento de Matem\'{a}ticas, Universidad
Aut\'{o}noma de Madrid, 28049 Madrid, Spain. (e-mail: fernando.quiros@uam.es).

\noindent\textsc{A. Rodr\'{\i}guez: } Departamento de Matem\'{a}tica Aplicada,  Universidad Polit\'{e}cnica de Madrid, 28040 Madrid, Spain.
(e-mail: ana.rodriguez@upm.es).

\end{document}